# Optimal Allocation of Series FACTS Devices Under High Penetration of Wind Power Within a Market Environment

Xiaohu Zhang, *Member, IEEE,* Di Shi, *Senior Member, IEEE,* Zhiwei Wang, *Member, IEEE,*
Bo Zeng, *Member, IEEE,* Xinan Wang, *Student Member, IEEE,* Kevin Tomsovic, *Fellow, IEEE,* Yanming Jin

*Abstract*—Series FACTS devices are one of the key enablers for very high penetration of renewables due to their capabilities in continuously controlling power flows on transmission lines. This paper proposes a bilevel optimization model to optimally locate variable series reactor (VSR) and phase shifting transformer (PST) in the transmission network considering high penetration of wind power. The upper level problem seeks to minimize the investment cost on series FACTS, the cost of wind power curtailment and possible load shedding. The lower level problems capture the market clearing under different operating scenarios. Due to the poor scalability of $B\theta$ formulation, the *shift factor* structure of FACTS allocation is derived. A customized reformulation and decomposition algorithm is designed and implemented to solve the proposed bilevel model with binary variables in both upper and lower levels. Detailed numerical results based on 118-bus system demonstrate the efficient performance of the proposed planning model and the important role of series FACTS for facilitating the integration of wind power.

## Nomenclature

*Indices*

| | |
|---|---|
| $i$, $j$ | Index of buses. |
| $k$ | Index of transmission elements. |
| $n$ | Index of generators. |
| $m$ | Index of loads. |
| $w$ | Index of wind farms. |
| $t$ | Index of scenarios. |

*Variables*

| | |
|---|---|
| $P_{nt}^g$ | Active power generation of generator $n$ in scenario $t$. |
| $P_{wt}^g$ | Active power generation of wind farm $w$ in scenario $t$. |

This project is funded by State Grid Corporation of China (SGCC) under project *Research on Spatial-Temporal Multidimensional Coordination of Energy Resources.*

Xiaohu Zhang, Di Shi, Zhiwei Wang are with GEIRINA, San Jose, CA, USA, email: {xiaohu.zhang,di.shi,zhiwei.wang}@geirina.net. Bo Zeng is with the Department of Industrial Engineering and the Department of Electrical & Computer Engineering, University of Pittsburgh, PA, USA, email: bzeng@pitt.edu. Xinan Wang is with the Department of Electrical Engineering, Southern Methodist University, Dallas, TX, USA, email: wx-napply@gmail.com. Kevin Tomsovic is with the Department of Electrical Engineering and Computer Science, the University of Tennessee, Knoxville, TN, USA, email: tomsovic@utk.edu. Yanming Jin is with State Grid Energy Research Institute, Beijing, China, email:jinyanming@sgeri.sgcc.cn.

| | |
|---|---|
| $P_{kt}$ | Active power flow on branch $k$ in scenario $t$. |
| $\Delta P_{mt}^d$ | Load shedding amount of load $m$ in scenario $t$. |
| $P_{wt}^{sp}$ | Wind power production spillage of wind farm $w$ in scenario $t$. |
| $\theta_k$ | The angle difference across branch $k$. |
| $\delta_k$ | Binary variable associated with placing a VSR on branch $k$. |
| $\alpha_k$ | Binary variable associated with placing a PST on branch $k$. |

*Parameters*

| | |
|---|---|
| $b_k$ | Negative susceptance of transmission line $k$. |
| $x_k$ | Reactance of transmission line $k$. |
| $a_n^g$ | Cost coefficient for generator $n$. |
| $\alpha$ | Cost coefficient for the wind power spillage. |
| $\beta$ | Cost coefficient for the load shedding. |
| $N_t$ | The number of operating hours of scenario $t$ in the target planning year. |
| $P_{nt}^{g,\min}$ | Minimum active power output of generator $n$ in scenario $t$. |
| $P_{nt}^{g,\max}$ | Maximum active power output of generator $n$ in scenario $t$. |
| $P_{mt}^d$ | Active power consumption of demand $m$ in scenario $t$. |
| $S_{kt}^{\max}$ | Capacity limit of branch $k$ in scenario $t$. |
| $P_{wt}^w$ | Available wind power production of wind farm $w$ in scenario $t$. |
| $C_k^V, A_k^V$ | Total and annualized investment cost of VSR on line $k$. |
| $C_k^P, A_k^P$ | Total and annualized investment cost of PST on line $k$. |
| $\boldsymbol{H}$ | Power transfer distribution factor matrix. |

*Sets*

| | |
|---|---|
| $\Omega_V, \Omega_P$ | Set of candidate transmission lines to install VSR and PST. |
| $\mathcal{D}$ | Set of loads. |
| $\mathcal{D}_i$ | Set of loads located at bus $i$. |
| $\Omega_L$ | Set of transmission lines. |
| $\Omega_T$ | Set of scenarios. |
| $\mathcal{B}$ | Set of buses. |
| $\mathcal{B}_{ref}$ | Set of reference bus. |
| $\mathcal{G}$ | Set of on-line generators. |

| $\mathcal{G}_i$ | Set of on-line generators located at bus $i$. |
|---|---|
| $\mathcal{W}$ | Set of wind farms. |
| $\mathcal{W}_i$ | Set of wind farms connected to bus $i$. |
| $f(i)$ | Set of lines specified as from bus $i$. |
| $t(i)$ | Set of lines specified as to bus $i$. |

Other symbols are defined as required in the text. Matrices are indicated by upper case bold, vectors by lower case bold.

# I. INTRODUCTION

WIND energy serves as one of the most effective approaches by power industry to reduce the emission of greenhouse gases to achieve sustainability. Many countries have set their own targets for wind power penetration levels in the future. For example, according to [1], the installed capacity of wind generation in the United States has reached over 61 GW in 2013 and it is anticipated that 20% and 30% of electricity demand will be supplied from wind power by 2030 and 2050, respectively.

There are many technical challenges related to the integration of wind power into the existing power grids. One of them is the intermittency of wind power generation, which brings more uncertainties to the power network and leads to potential risks on the system reliability and efficiency [2]. Another issue is the insufficient transmission capacity to deliver large amount of wind power from remote areas to the load centers [3]. An obvious approach to enhance the transmission capacity and facilitate the wind power integration is to construct new transmission lines. Nevertheless, such projects are unattractive since they usually require high investment cost, long construction time and stringent environmental approvals [4]. An economical approach is to install Flexible AC Transmission Systems (FACTS) on selected lines, which improves the utilization of the existing transmission infrastructure by regulating power flows in a more flexible way [5]. With the power flow control capability introduced by FACTS, transmission bottlenecks can be avoided through shifting the power from the congested lines to the underutilized lines nearby. In addition, due to their fast operations, i.e., often within a few cycles of system frequency, FACTS can be dynamically adjusted to accommodate the stochastic nature of wind power [3]. Further, it is expected that more FACTS-like devices [6] with much lower price will be commercially available soon under the efforts of Green Electricity Network Integration (GENI) program [7]. Thus, efficient planning models, which provide useful information regarding the optimal locations of FACTS devices, should be developed to facilitate the integration of increasing wind power.

In the technical literature, determining the optimal locations and compensation levels of FACTS devices have been studied extensively. Since the mathematical formulations are originally nonlinear and non-convex, various heuristic methods have been proposed to solve the FACTS allocation problem. The authors in [4] leverage the genetic algorithm (GA) to determine the optimal placements of phase shifting transformer (PST) and thyristor controlled series compensator (TCSC). In [8], a particle swarm optimization (PSO) is proposed to identify the locations of TCSC to enhance the loadability of the system. In [9], a hybrid PSO/SQP algorithm is used to find the optimal placements of TCSCs and their compensation levels under different operating conditions accordingly. Reference [10] provides a review regarding heuristic methods application in FACTS allocation problem. Priority indices approaches [11], [12] are another category of FACTS placement methods. These approaches derive certain types of indices which indicate the impacts of FACTS devices on different system objectives, such as, loss minimization, congestion relief and transfer capability enhancement, and so on.

Recently, mixed integer program (MIP) based approaches have also been proposed to allocate FACTS. The authors in [13] formulate the PST placement problem as a mixed integer linear program (MILP) to improve the system loadability. Reference [14] leverages the line flow based equations to allocate TCSC via MILP. The nonlinear term in the original model is approximately linearized by relaxing one variable to its hard limits. Reference [15] proposes an MILP based planning model to place TCSC considering multi-scenarios including both base operating states and contingencies. The nonlinear part introduced by the variable reactance is exactly linearized by a reformulation technique. In [16], the authors propose a two-stage stochastic programming to co-optimize the locations of TCSC and transmission switch considering wind power uncertainty. The model is solved by branch-and-price algorithm.

It is well understood that the investment made by the system planner should facilitate the energy trading in the electricity market. To explicitly incorporate the market clearing conditions in the planning model, bilevel optimization is usually used. Previous studies have leveraged the bilevel model in various types of power system planning problem such as transmission expansion planning (TEP) [17], [18], wind farm investment [19] and conventional generator planning [20]. In all these aforementioned models, the market clearing conditions under different scenarios are represented by a collection of lower level problems.

In this paper, we propose a bilevel model to co-optimize the locations of two types of series FACTS: variable series reactor (VSR) and PST. These two devices can efficiently vary reactance and phase angle difference across the transmission line respectively, so as to regulate the power flow. They are suitable for power system congestion relief and wind power integration. To achieve a desired compromise between computational tractability and model accuracy, we adopt a static model which focuses on a single representative year [21]. Fig. 1 provides the general structure of the proposed model. In our bilevel model, the system planner in the upper level aims at finding the optimal locations of series FACTS to minimize the wind power curtailment, involuntary load shedding and annualized investment cost in the single target year, subject to maximum number of devices for each type of series FACTS. The amount of wind power curtailment and load shedding are determined by a series of lower level problems representing the market clearing for different load-wind scenarios, with the consideration of FACTS operations. In addition, due to the poor scalability of the $B\theta$ formulation, commercial solvers used by power industry usually adopt the *shift factor* formulation using power transfer distribution factor (PTDF) [22]–[24]. Hence, the *shift factor* structure for the allocation of series FACTS is also derived.

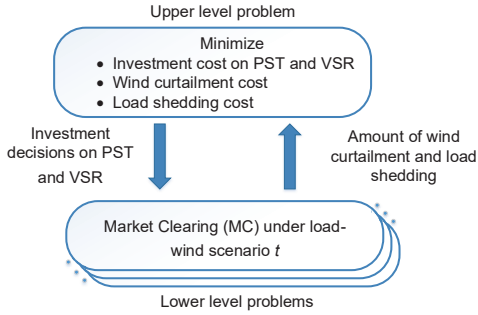

Fig. 1. Bilevel structure of proposed series FACTS investment model.

Note that the market clearing in each lower level problem usually adopts DC optimal power flow (OPF) model so it is continuous and linear. Hence, the complete bilevel model can be reformulated as a single level problem by replacing each lower problem with its Karush-Kuhn-Tucker (KKT) conditions or primal-dual formulation. Nevertheless, the lower level problems in our proposed model are non-convex and nonlinear due to the variable susceptance by VSR. We first leverage an exact reformulation technique to linearize the nonlinear part, which introduces additional dummy binary variables in the lower level problems. Then a recently proposed reformulation and decomposition algorithm [25] is adaptively developed to handle this challenging MIP bilevel model with binary variables in both levels.

Considering the extensive literature in this area, the main contributions of this paper are:

1) a stochastic MIP bilevel model to co-optimize the locations of VSR and PST using *shift factor* structure within a market environment under very high levels of wind power;

2) a decomposition algorithm to solve the stochastic MIP bilevel model;

3) detailed results with respect to the significant benefits of series FACTS in reducing the amount wind curtailment; and

4) analysis of the computational performance of $B\theta$ and *shift factor* formulation in terms of solution time and model size for the proposed FACTS allocation model.

The rest of this paper is organized as follows. In section II, the injection models of series FACTS and the reformulation technique are presented. Section III illustrates the formulation of the proposed bilevel model. The solution approach based on the decomposition algorithm is demonstrated in Section IV. Section V provides the numerical results based on 118-bus system. Finally, conclusions are given in section VI.

## II. Power Injection models of Series FACTS

The PTDF is defined as the sensitivity of the power flow on line $k$ with respect to the power injection at bus $i$, which is often computed offline with system topology and branch impedance [23]. Obviously, the installation of FACTS can change the susceptance of the transmission line and lead to variable PTDF. To make the PTDF constant, the series FACTS are modeled by power injections at both ends of the selected transmission lines [24].

### A. Injection model of PST

In steady state, the classic model of PST on line $k$ is represented by a continuously variable phase angle $\theta_k^P$ in series with the line reactance $x_k$. Fig. 2 depicts the transformation from the PST classic model to its power injection model. Mathematically, this conversion can be expressed as [26]:

$$\tilde{P}_k = b_k(\theta_k + \alpha_k \theta_k^P)$$
$$= b_k\theta_k + \alpha_k b_k \theta_k^P = P_k + \psi_k^P \qquad (1)$$

where $\tilde{P}_k$ and $P_k$ are the power flows on line $k$ with and without PST; $\psi_k^P$ is the active power injection introduced by PST. Note that $\psi_k^P$ involves the product between a binary variable and a continuous variable, which can be linearized by (2):

$$\alpha_k b_k \theta_k^{P,\min} \le \psi_k^P \le \alpha_k b_k \theta_k^{P,\max} \qquad (2)$$

Specifically, when line $k$ is not selected to install PST, i.e., $\alpha_k = 0$, the active power injection $\psi_k^P$ will be zero; otherwise, the power injection $\psi_k^P$ will be bounded by its upper and lower limits.

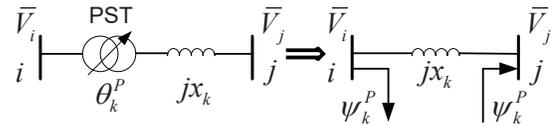

Fig. 2. Conversion of PST classic model to its power injection model.

### B. Injection model of VSR

Similarly, Fig. 3 demonstrates the transformation from VSR classic model, i.e., variable reactance $x_k^V$ in series with $x_k$, to its power injection model.

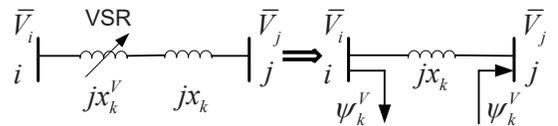

Fig. 3. Conversion of VSR classic model to its power injection model.

This conversion can be mathematically derived as

$$\tilde{P}_k = (b_k + \delta_k b_k^V)\theta_k$$
$$= b_k\theta_k + \delta_k b_k^V \theta_k = P_k + \delta_k \frac{b_k^V}{b_k}P_k = P_k + \psi_k^V \qquad (3)$$

where $\tilde{P}_k$ and $P_k$ are the power flows on line $k$ with and without VSR; $b_k^V$ is the susceptance change introduced by VSR; $\psi_k^V$ denotes the power injection given by VSR. Note that $\psi_k^V$ is a trilinear term because it is the product of one binary variable and two continuous variables. To linearize this term, we first define $\frac{b_k^V}{b_k} = \Delta b_k^V$ and the power injection $\psi_k^V$ can be written as:

$$\psi_k^V = \delta_k \Delta b_k^V P_k \tag{4}$$

With the upper and lower bounds of $\Delta b_k^V$, constraint (4) can be expressed as

$$\delta_k \Delta b_k^{V,\min} \leq \frac{\psi_k^V}{P_k} = \delta_k \Delta b_k^V \leq \delta_k \Delta b_k^{V,\max} \tag{5}$$

Since the sign of $P_k$ cannot be determined beforehand, we introduce a binary variable $u_k$ and leverage the big-M complementary constraints [27] to rewrite (5) into (6) and (7):

$$- M_{k1} u_k + \delta_k P_k \Delta b_k^{V,\min} \leq \psi_k^V \leq \delta_k P_k \Delta b_k^{V,\max} + M_{k1} u_k \tag{6}$$

$$- M_{k1}(1 - u_k) + \delta_k P_k \Delta b_k^{V,\max} \leq \psi_k^V$$
$$\leq \delta_k P_k \Delta b_k^{V,\min} + M_{k1}(1 - u_k) \tag{7}$$

Due to the sufficiently large number $M_{k1}$, only one constraint of (6) and (7) will be active during the optimization process and the other one will become a redundant one. It should be noted that the term $\delta_k P_k$ in (6) and (7) is nonlinear. We introduce another continuous variable $v_k = \delta_k P_k$ and linearize it by using the big-M method:

$$- \delta_k M_{k2} \leq v_k \leq \delta_k M_{k2} \tag{8}$$

$$P_k - (1 - \delta_k) M_{k2} \leq v_k \leq P_k + (1 - \delta_k) M_{k2} \tag{9}$$

Then constraint (6) and (7) can be written as:

$$- M_{k1} u_k + v_k \Delta b_k^{V,\min} \leq \psi_k^V \leq v_k \Delta b_k^{V,\max} + M_{k1} u_k \tag{10}$$

$$- M_{k1}(1 - u_k) + v_k \Delta b_k^{V,\max} \leq \psi_k^V$$
$$\leq v_k \Delta b_k^{V,\min} + M_{k1}(1 - u_k) \tag{11}$$

Hence, the nonlinear power injection $\psi_k^V$ is linearized by using (8)-(11).

## III. Problem Formulation

The bilevel model for series FACTS investment is comprised by an upper level problem and a collection of lower level problems. The upper level problem seeks to minimize the wind curtailment, involuntary load shedding and investment cost in series FACTS for a single target year. With the investment decisions from the upper level problem, each of the lower level problem, one per load-wind scenario, represents the market clearing conditions using *shift factor* formulation. This bilevel model is presented as follows:

$$\min_{\boldsymbol{x} \cup \boldsymbol{y} \cup \boldsymbol{z}} \sum_{k \in \Omega_V} A_k^V \delta_k + \sum_{k \in \Omega_V} A_k^P \alpha_k$$
$$+ \alpha \sum_{t \in \Omega_T} N_t \sum_{w \in \mathcal{W}} P_{wt}^{sp} + \beta \sum_{t \in \Omega_T} N_t \sum_{m \in \mathcal{D}} \Delta P_{mt}^d \tag{12a}$$

subject to

$$\sum_{k \in \Omega_V} \delta_k \leq N_V \tag{12b}$$

$$\sum_{k \in \Omega_P} \alpha_k \leq N_P \tag{12c}$$

$$\delta_k + \alpha_k \leq 1, \ \forall k \in \Omega_V \cap \Omega_P \tag{12d}$$

where $P_{wt}^{sp}$ and $\Delta P_{mt}^d \in \arg\{$

$$\min_{\boldsymbol{y} \cup \boldsymbol{z}} \sum_{n \in \mathcal{G}} a_n^g P_{nt}^g + \beta \sum_{m \in \mathcal{D}} \Delta P_{mt}^d \tag{12e}$$

s.t $(2), (8) - (11)$ and

$$P_{kt} = \boldsymbol{H}(k,i)(\sum_{n \in \mathcal{G}_i} P_{nt}^g + \sum_{w \in \mathcal{W}_i} P_{wt}^g - \sum_{m \in \mathcal{D}_i} (P_{mt}^d - \Delta P_{mt}^d)$$
$$- \sum_{k \in \Omega_P^{i(fr)}} \psi_{kt}^P + \sum_{k \in \Omega_P^{i(to)}} \psi_{kt}^P - \sum_{k \in \Omega_V^{i(fr)}} \psi_{kt}^V + \sum_{k \in \Omega_V^{i(to)}} \psi_{kt}^V),$$
$$\forall k, \forall t \tag{12f}$$

$$\sum_{n \in \mathcal{G}} P_{nt}^g + \sum_{w \in \mathcal{W}} P_{wt}^g - \sum_{m \in \mathcal{D}} (P_{mt}^d - \Delta P_{mt}^d) = 0, \ \forall t \tag{12g}$$

$$P_{wt}^{sp} = P_{wt}^a - P_{wt}^g, \ \forall w, \forall t \tag{12h}$$

$$0 \leq P_{wt}^g \leq P_{wt}^a, \ \forall w, \forall t \tag{12i}$$

$$P_{nt}^{g,\min} \leq P_{nt}^g \leq P_{nt}^{g,\max}, \ \forall n, \forall t \tag{12j}$$

$$- S_{kt}^{\max} \leq P_{kt} \leq S_{kt}^{\max}, \ \forall k \setminus k \in (\Omega_V \cup \Omega_P), \forall t \tag{12k}$$

$$- S_{kt}^{\max} \leq P_{kt} + \psi_{kt}^V \leq S_{kt}^{\max}, \ \forall k \in \Omega_V, \forall t \tag{12l}$$

$$- S_{kt}^{\max} \leq P_{kt} + \psi_{kt}^P \leq S_{kt}^{\max}, \ \forall k \in \Omega_P, \forall t \tag{12m}$$

$$0 \leq \Delta P_{mt}^d \leq P_{mt}^d, \ \forall m, \forall t \ \} \tag{12n}$$

The continuous optimization variables of the lower level problems comprise the elements in set $\boldsymbol{y} = \{P_{nt}^g, P_{wt}^g, P_{wt}^{sp}, \Delta P_{mct}^d, P_{kt}, \psi_{kt}^V, \psi_{kt}^P, v_{kt}\}$. The binary variables of the lower level problems are in set $\boldsymbol{z} = \{u_{kt}\}$. The upper level decision variables are represented by $\boldsymbol{x} = \{\alpha_k, \delta_k\}$.

The focus of this paper is to minimize the wind power spillage and to determine whether the wind power should be curtailed, or series FACTS should be installed from the economic point of view. Moreover, the reliability issue is one of the primary concerns for the system planner. Thus, in the upper level problem, the objective function (12a) seeks to minimize the annualized investment cost in VSR and PST (first two terms) plus the annual wind curtailment (third term) and load shedding cost (fourth term). The amount of wind curtailment and load shedding in the target planning year are computed by multiplying $P_{wt}^{sp}$ and $\Delta P_{mt}^d$ in each scenario $t$ with their corresponding operating hours $N_t$. Constraint (12b) and (12c) limit the number of VSRs and PSTs that can be installed in the system, respectively. Constraint (12d) denotes that a line can be equipped with a VSR or PST but not both.

The upper level problem is also constrained by a collection of lower level problems which represent the market clearing conditions under different load-wind scenarios. For each scenario $t$, the objective function (12e) is to minimize the production cost from the conventional generators and the possible load shedding. This is equivalent to maximize the social welfare if the demand is considered inelastic [21]. Constraint (12f) denotes the power flow formulation using the PTDF matrix $\boldsymbol{H}$. Note that the installation of series FACTS injects active power at one end of the selected transmission line and withdraws power at the other end. So the total power injection/consumption at bus $i$ should be modified to include the possible power injection/consumption from series FACTS. The sets $\Omega_P^{i(fr)}$ and $\Omega_P^{i(to)}$ in (12f) represent the candidate lines to install PST with their from and to buses to be $i$. Similar

 

description applies to $\Omega_V^{i(fr)}$ and $\Omega_V^{i(to)}$. The power balance constraint is enforced by (12g). Constraint (12h) represents that the wind power spillage is computed as the difference between the available wind power and the dispatched wind power. Constraints (12i) states that the dispatched wind power is bounded by its available amount. The generation limits of generators are enforced by constraint (12j). Constraint (12k)-(12m) consider the thermal limits of normal lines, candidate lines to install VSR and PST, respectively. Note that by following [19]–[21], we adopt thermal limits to bound the power flow on all lines for simplicity. The power flow on medium or long transmission lines can be bounded by voltage or angular stability limits. Note also that the flow limit of a line compensated by VSR can increase from voltage or angular stability limit to thermal limit [28]. Finally, constraint (12n) enforces an upper bound on the load curtailment.

## IV. Solution Approach

As mentioned in the introduction, the proposed bilevel model cannot be directly solved by the KKT or primal-dual based method due to the existence of the binary variables in the lower level problems. We then leverage a recently proposed reformulation and decomposition algorithm to address this challenge. For simplicity, our proposed bilevel model is first compactly written as:

$$\min_{\boldsymbol{x} \cup \boldsymbol{y} \cup \boldsymbol{z}} \quad \boldsymbol{f}^T \boldsymbol{x} + \sum_{t \in \Omega_T} N_t \boldsymbol{g}_t^T \boldsymbol{y}_t \tag{13a}$$

$$\boldsymbol{A}\boldsymbol{x} \leq \boldsymbol{b} \tag{13b}$$

where $\boldsymbol{y}_t$ and $\boldsymbol{z}_t \in \arg\{ \min_{\boldsymbol{y}_t \cup \boldsymbol{z}_t} \quad \boldsymbol{w}_t^T \boldsymbol{y}_t \tag{13c}$

$$\text{s.t.} \quad \boldsymbol{E}\boldsymbol{y}_t = \boldsymbol{h}_t \tag{13d}$$

$$\boldsymbol{P}\boldsymbol{y}_t + \boldsymbol{Q}\boldsymbol{z}_t \leq \boldsymbol{r}_t - \boldsymbol{K}\boldsymbol{x} \} \quad \forall t \in \Omega_T \tag{13e}$$

In (13), $\boldsymbol{x}$ represents the upper level decision variables, $\boldsymbol{y}_t$ and $\boldsymbol{z}_t$ are defined as the continuous and binary variables in the lower level problem under scenario $t$.

Next, the following equivalent formulation of (13) is obtained by duplicating the decision variables and constraints in the lower level problems [25]:

$$\min_{\boldsymbol{x} \cup \boldsymbol{y} \cup \boldsymbol{z}} \quad \boldsymbol{f}^T \boldsymbol{x} + \sum_{t \in \Omega_T} N_t \boldsymbol{g}_t^T \tilde{\boldsymbol{y}}_t \tag{14a}$$

$$\boldsymbol{A}\boldsymbol{x} \leq \boldsymbol{b} \tag{14b}$$

$$\boldsymbol{E}\tilde{\boldsymbol{y}}_t = \boldsymbol{h}_t \tag{14c}$$

$$\boldsymbol{P}\tilde{\boldsymbol{y}}_t + \boldsymbol{Q}\tilde{\boldsymbol{z}}_t \leq \boldsymbol{r}_t - \boldsymbol{K}\boldsymbol{x} \tag{14d}$$

$$\boldsymbol{w}_t^T \tilde{\boldsymbol{y}}_t \leq \min_{\boldsymbol{y}_t \cup \boldsymbol{z}_t} \{ \quad \boldsymbol{w}_t^T \boldsymbol{y}_t \tag{14e}$$

$$\text{s.t.} \quad \boldsymbol{E}\boldsymbol{y}_t = \boldsymbol{h}_t \tag{14f}$$

$$\boldsymbol{P}\boldsymbol{y}_t + \boldsymbol{Q}\boldsymbol{z}_t \leq \boldsymbol{r}_t - \boldsymbol{K}\boldsymbol{x} \} \quad \forall t \in \Omega_T \tag{14g}$$

We mention that constraint (14c)-(14g) guarantee that, for a given $\boldsymbol{x}$, $(\tilde{\boldsymbol{y}}_t, \tilde{\boldsymbol{z}}_t)$ is not only feasible but also optimal to the lower-level problem [25]. Hence, (13) and (14) are equivalent. Noting that $\boldsymbol{z}_t$ is in a finite binary set for any $t$, we can rewrite (14e)-(14g) by enumerating all possible values of $\boldsymbol{z}_t$ as following.

$$\boldsymbol{w}_t^T \tilde{\boldsymbol{y}}_t \leq \min_{\boldsymbol{y}_t} \{ \quad \boldsymbol{w}_t^T \boldsymbol{y}_t : \tag{15a}$$

$$\text{s.t.} \quad \boldsymbol{P}\boldsymbol{y}_t \leq \boldsymbol{r}_t - \boldsymbol{K}\boldsymbol{x} - \boldsymbol{N}\boldsymbol{z}_t^{l^*} \tag{15b}$$

$$\boldsymbol{E}\boldsymbol{y}_t = \boldsymbol{h}_t \quad \} \quad \forall \boldsymbol{z}_t^{l^*} \in \mathcal{Z}_t \tag{15c}$$

where $\mathcal{Z}_t$ is the collection of all possible $\boldsymbol{z}_t$ and $\boldsymbol{z}_t^{l^*}$ is a particular realization of $\boldsymbol{z}_t$. Two observations can be drawn from (15). First, given a fixed $\boldsymbol{z}_t^{l^*}$, the optimization problem in the "min" operator of (15) becomes a linear program (LP), which can be simply and equivalently replaced by primal-dual or KKT conditions. Second, a partial enumeration leads to a relaxation of (14). So, starting from an initial value of $\boldsymbol{z}_t$, we progressively add more realizations of $\boldsymbol{z}_t$ into (14) and solve tighter relaxations of (14). Overall, a decomposition approach based on the widely adopted column-and-constraint generation algorithm [29]–[31] can be developed.

### A. Decomposition Algorithm: Column-and-Constraint Generation Method

The column-and-constraint generation based decomposition algorithm involves solving one master problem and two subproblems iteratively. Fig. 4 depicts the flowchart of the algorithm.

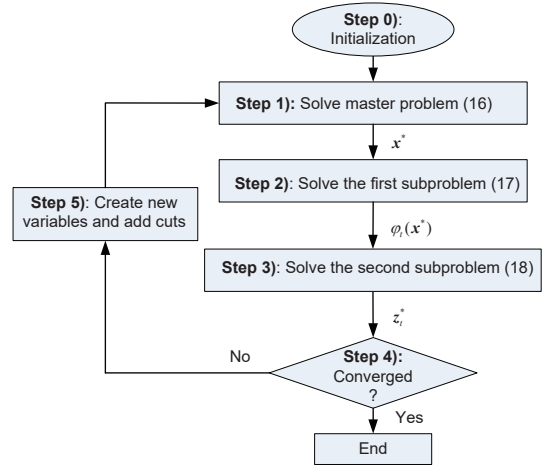

Fig. 4. Flowchart of the decomposition algorithm. Details of each step are illustrated in Section IV-A.

The complete procedures are described as follows:

0) Set $LB = -\infty$, $UB = +\infty$ and iteration counter $q = 0$. Select a small tolerance $\epsilon$ to control the convergence.

1) Solve the master problem (**MP**) $\forall t \in \Omega_T$:

$$\Phi = \min_{\boldsymbol{x} \cup \boldsymbol{y} \cup \boldsymbol{z}} \quad \boldsymbol{f}^T \boldsymbol{x} + \sum_{t \in \Omega_T} N_t \boldsymbol{g}_t^T \tilde{\boldsymbol{y}}_t \tag{16a}$$

$$\text{s.t.} \quad (14b) - (14d) \text{ and}$$

$$\boldsymbol{w}_t^T \tilde{\boldsymbol{y}}_t \leq (\boldsymbol{z}_t^{l^*})^T \boldsymbol{Q}^T - \boldsymbol{r}_t^T) \boldsymbol{\lambda}_t^{(l)} - \boldsymbol{h}_t^T \boldsymbol{\mu}_t^{(l)}$$
$$+ \boldsymbol{x}^T \boldsymbol{K}^T \boldsymbol{\lambda}_t^{(l)}, \ 1 \leq l \leq q \tag{16b}$$

$$\boldsymbol{P}^T \boldsymbol{\lambda}_t^{(l)} + \boldsymbol{E}^T \boldsymbol{\mu}_t^{(l)} + \boldsymbol{w}_t = 0, \ 1 \leq l \leq q \tag{16c}$$

$$\boldsymbol{\lambda}_t^{(l)} \geq 0, \ 1 \leq l \leq q \tag{16d}$$

Derive its optimal solution, obtain the values of upper level decision variables $\boldsymbol{x}^*$ and update $LB = \Phi$.

2) With $\boldsymbol{x}^*$ from step 1), solve the first subproblem (**SP1**) for every $t$ in $\Omega_T$.

$$\varphi_t(\boldsymbol{x}^*) = \min_{\boldsymbol{y}_t \cup \boldsymbol{z}_t} \quad \boldsymbol{w}_t^T \boldsymbol{y}_t \tag{17a}$$

$$\text{s.t.} \quad \boldsymbol{E}\boldsymbol{y}_t = \boldsymbol{h}_t \tag{17b}$$

$$\boldsymbol{P}\boldsymbol{y}_t + \boldsymbol{Q}\boldsymbol{z}_t \le \boldsymbol{r}_t - \boldsymbol{K}\boldsymbol{x}^* \tag{17c}$$

Derive their optimal values $\varphi_t(\boldsymbol{x}^*)$.

3) With $\boldsymbol{x}^*$ and $\varphi_t(\boldsymbol{x}^*)$, solve the second subproblem (**SP2**) $\forall t \in \Omega_T$.

$$\phi(\boldsymbol{x}^*) = \min_{\boldsymbol{y} \cup \boldsymbol{z}} \quad \sum_{t \in \Omega_T} N_t \boldsymbol{g}_t^T \boldsymbol{y}_t \tag{18a}$$

$$\text{s.t.} \quad \boldsymbol{w}_t^T \boldsymbol{y}_t \le \varphi_t(\boldsymbol{x}^*) \tag{18b}$$

$$\boldsymbol{E}\boldsymbol{y}_t = \boldsymbol{h}_t \tag{18c}$$

$$\boldsymbol{P}\boldsymbol{y}_t + \boldsymbol{Q}\boldsymbol{z}_t \le \boldsymbol{r}_t - \boldsymbol{K}\boldsymbol{x}^* \tag{18d}$$

Derive its optimal solution $(\boldsymbol{y}_t^*, \boldsymbol{z}_t^*)$ and update $UB = \min\{UB, \boldsymbol{f}^T \boldsymbol{x}^* + \phi(\boldsymbol{x}^*)\}$.

4) If $\left|\frac{UB - LB}{UB}\right| \le \epsilon$, return $UB$ and the corresponding solutions. Stop the algorithm. Otherwise, go to step 5).

5) Set $\boldsymbol{z}_t^{(q+1)} = \boldsymbol{z}_t^*$, create new variables $\boldsymbol{\mu}_t^{(q+1)}, \boldsymbol{\lambda}_t^{(q+1)}$ and add the following constraints (cuts) to **MP**:

$$\boldsymbol{w}_t^T \tilde{\boldsymbol{y}}_t \le (\boldsymbol{z}_t^{(q+1)})^T \boldsymbol{Q}^T - \boldsymbol{r}_t^T) \boldsymbol{\lambda}_t^{(q+1)} - \boldsymbol{h}_t^T \boldsymbol{\mu}_t^{(q+1)} + \boldsymbol{x}^T \boldsymbol{K}^T \boldsymbol{\lambda}_t^{(q+1)} \tag{19a}$$

$$\boldsymbol{P}^T \boldsymbol{\lambda}_t^{(q+1)} + \boldsymbol{E}^T \boldsymbol{\mu}_t^{(q+1)} + \boldsymbol{w}_t = 0 \tag{19b}$$

$$\boldsymbol{\lambda}_t^{(q+1)} \ge 0 \tag{19c}$$

Set $q = q + 1$ and go to step 1).

Note that $\boldsymbol{\mu}_t$ and $\boldsymbol{\lambda}_t$ are the dual variables associated with the equality and inequality constraints in the lower level problems. With a given $\boldsymbol{z}_t$, we replace each of the lower level problem with its primal-dual reformulation because it is computationally more friendly than the KKT based one [25]. In (19a), there is a nonlinear term that is a product of $\boldsymbol{\lambda}_t$ and $\boldsymbol{x}$. Fortunately, $\boldsymbol{x}$ only comprises binary variables, i.e., investment decisions on series FACTS, so this nonlinear term can also be linearized by using big-M method.

### B. Computational Enhancement

*1) Candidate Locations Selection:* In real power system, considering every transmission line as a candidate location for series FACTS is impractical and unnecessary. For this reason, we first perform a preliminary experiment based on the sensitivity approach [11] to obtain the candidate lists of VSR and PST.

The procedures to determine the VSR candidate locations are provided below:

1) Run a DCOPF for each load-wind scenario without VSR. Every line reactance is treated as an optimization variable and the following constraint is included in the OPF model, i.e., fix the line reactance to its original value:

$$\tilde{x}_{kt} = x_k \tag{20}$$

Note that the OPF model is nonlinear so IPOPT [32] is leveraged to solve it.

2) Obtain the sensitivity ($\eta_{kt}$) of the operation cost with respect to the change of line reactance in each scenario, i.e., value of the dual variable associated with constraint (20).

3) Compute the weighted sensitivity ($\bar{\eta}_k$) of branch $k$ by equation (21):

$$\bar{\eta}_k = \sum_{t \in \Omega_T} N_t |\eta_{kt} x_k| \tag{21}$$

4) Sort $\bar{\eta}_k$ in a descending order and select the first 10 lines as candidate locations for VSR.

Similar procedures are applied to obtain the candidate locations for PST.

*2) Lower Level Problems Size Reduction:* Given $\mathcal{Z}_t$ a finite binary set, the decomposition algorithm converges in finite iterations [25]. The number of binary variables in the lower level problems has a large impact on the computational burden. From (6) and (7), the binary variable $u_k$ indicates the sign of $P_k$, which is the power flow direction of the VSR candidate line $k$. Based on the engineering insights, most of the lines to be equipped with series FACTS are usually tie lines whose flow directions are not likely to vary [24], [33]. For instance, the flow direction on the California Oregon Intertie (COI) can be easily predicted.

In addition, one advantage of *shift factor* formulation over $B\theta$ formulation is that the *shift factor* structure allows the system operator to monitor a subset of transmission lines that are interested, e.g., lines over a certain voltage level, lines that are usually overloaded based on the historical data, etc. From the modeling point of view, this can be achieved by selecting the corresponding rows of the PTDF matrix $\boldsymbol{H}$. In power industry, this feature has already been implemented by several commercial planning softwares. As an example, PLEXOS has a user option to exclude a set of transmission lines in the long term planning process [34].

Hence, based on the above two observations, we propose the following procedures to reduce the size of lower level problems:

1) Run DCOPF without series FACTS for each load-wind scenario.

2) Identify the power flow directions of the VSR candidate lines in each scenario.

3) Obtain the thermal loadings of normal lines, i.e., not the candidate lines, in each scenario.

4) For the VSR candidate lines whose flow directions do not change among all the scenarios, we fix their flow directions to the OPF results in the lower level problems.

5) For those normal lines whose thermal loadings are below 60% in all the scenarios, we exclude them in our planning model.

## V. NUMERICAL RESULTS

The proposed model and solution approach are tested on the IEEE 118-bus system. The system data can be found in [35]. As one typical example of VSR, TCSC is selected in the case studies. It is assumed that the compensation range of TCSC varies from -70% to 20% of its corresponding line reactance [9]. Moreover, we assume the range of phase shift angle to be $[-10°, +10°]$. According to [26], the investment cost of a PST is dependent on rating of the line concerned. The cost coefficient is selected to be 100\$/kVA [36] so the total investment cost of a PST, i.e., $C_k^P$, can be expressed as:

$$C_k^P = 100 \cdot S_k^{\max} \cdot 1000 \quad (22)$$

In terms of \$/kVar, the investment of a TCSC can be achieved by [9]:

$$I_k^V = 0.0015(S_k^V)^2 - 0.713 S_k^V + 153.75 \quad (23)$$

where $I_k^V$ is the cost in \$/kVar and $S_k^V$ is the maximum compensation level of the device in MVar, which can be expressed as [37], [38]:

$$S_k^V = \frac{(S_k^{\max})^2}{S_b} \overline{x}_k^V \quad (24)$$

where $S_b$ is the base MVA for the system, $\overline{x}_k^V$ is the maximum reactance that a TCSC can compensate. Thus, the investment cost for a TCSC can be written as:

$$C_k^V = I_k^V \cdot S_k^V \cdot 1000 \quad (25)$$

Note that the annualized investment cost of FACTS device is computed by its total cost along with the interest rate and life time of the device by using the following equations [39]:

$$A_k^P = C_k^P \cdot \frac{d(1+d)^{LT}}{(1+d)^{LT}-1} \quad (26)$$

$$A_k^V = C_k^V \cdot \frac{d(1+d)^{LT}}{(1+d)^{LT}-1} \quad (27)$$

where $d$ is the yearly interest rate and $LT$ is the life time of the device. In this work, $LT$ is selected to be 5 years and $d$ is 5% [16], [39]. The cost coefficient of wind curtailment $\alpha$ is assumed to be 50 \$/MWh [18]. The load shedding penalty coefficient $\beta$ is set to be 5000 \$/MWh. Finally, through a trial and error process, $M_{k1}$ and $M_{k2}$ are selected to be $2 \cdot \max(|\Delta b_k^{V,\min}|, |\Delta b_k^{V,\max}|) \cdot S_k^{\max}$ and $3.5 S_k^{\max}$, respectively. The calculation of the range of $\Delta b_k^V$ is provided in Appendix A.

To obtain the load-wind scenarios, we assume that the annual load of the test system follows the normalized load profile from the 2015 ISO New England hourly demand reports [40]. In addition, the hourly wind power intensities provided by [41] are used to represent the wind generation profile. We then use *K-means* method [42] to conduct the scenario reduction. Nevertheless, the extreme operating conditions may be eliminated using *K-means* method by selecting cluster centroids. Before implementing *K-means*, we first extract two extreme conditions corresponding to the highest demand and highest wind generation levels. Then *K-means* clustering is leveraged to reduce the number of scenarios from 8758 to 18. The number of operating hours, load levels and wind intensities for the final 20 scenarios are provided in Table I.

| # | # of hours | Load levels | Wind intensities | # | # of hours | Load levels | Wind intensities |
|---|---|---|---|---|---|---|---|
| 1 | 486 | 0.4858 | 0.3023 | 11 | 356 | 0.5323 | 0.7927 |
| 2 | 391 | 0.6916 | 0.8007 | 12 | 202 | 0.8558 | 0.1858 |
| 3 | 361 | 0.7338 | 0.6263 | 13 | 677 | 0.6266 | 0.5018 |
| 4 | 690 | 0.5919 | 0.0825 | 14 | 463 | 0.4948 | 0.4203 |
| 5 | 561 | 0.4796 | 0.1846 | 15 | 120 | 0.9065 | 0.5088 |
| 6 | 452 | 0.4870 | 0.5815 | 16 | 423 | 0.7437 | 0.4031 |
| 7 | 410 | 0.7026 | 0.2600 | 17 | 822 | 0.5897 | 0.2117 |
| 8 | 760 | 0.6036 | 0.3488 | 18 | 440 | 0.7087 | 0.1213 |
| 9 | 503 | 0.4701 | 0.0844 | 19 | 1 | 1.0000 | 0.1840 |
| 10 | 641 | 0.5936 | 0.6530 | 20 | 1 | 0.4915 | 0.8670 |

### A. IEEE 118-Bus System

The IEEE 118-bus system has 19 generators and 185 transmission lines. The peak loads are assumed to be 1.2 times their values provided in [35]. In addition, the thermal limits for the transmission lines decrease to 75% of the values in [35]. It is assumed that three wind farms, with the maximum capacity of 1600 MW, are located at bus 5, 26 and 91 [16], [43]. We also assume that the wind power intensities of wind farm at bus 5 and 26 are the values provided by Table I and wind intensities of wind farm at bus 91 are 10% lower than the values given in Table I. The number of candidate lines to install TCSC and PST are both selected to be 10. Among the TCSC candidate lines, the flow directions of eight lines are fixed so the number of binary variables in the lower level problems is reduced from 200 to 40.

Table II provides the planning results for seven cases regarding various limits on the number of TCSCs and PSTs. Column 2-4 represent the annual amount of wind curtailment for each wind farm. The fifth column shows the annual load shedding amount. The locations and annualized investment cost of TCSC and PST are given in column 6-9. Column 10 provides the value of objective function, i.e., (12a). The eleventh column indicates the wind penetration level, which is defined as the portion of load that can be covered by wind generation on an annual basis. The last column gives the computational time.

As observed from Table II, without any series FACTS, the amount of wind curtailment in the target planning year is 5.56e6 MWh and this value decreases to 4.80e6 MWh and 4.53e6 MWh with one TCSC and PST respectively. Although the PST has a higher investment cost than a TCSC, it has more impacts on the integration of wind power. In addition, the installation of one TCSC or PST can both eliminate the load shedding. As the maximum number of TCSC ($N_V$) and PST ($N_P$) increase, the objective value decreases and the wind penetration level increases. When comparing the value of objective function for the case without series FACTS and the case with two TCSCs and PSTs, a total savings of \$84.78 M can be achieved. Moreover, the wind penetration level is increased by 6.80%. The result indicates that 6.80% of the total loads which was served by conventional generators can be provided by cheap wind generations.

Note that if constraint (12b) and (12c) are eliminated, i.e., no limits on the number of series FACTS, the planning model suggests installing seven TCSCs and three PSTs. The value of objective function is $174.28 M with a wind penetration level of 42.03%. The amounts of wind curtailment at bus 5, 26 and 91 are 0.82e6 MWh, 1.84e6 MWh and 0.44e6 MWh, respectively. If more FACTS devices are allowed to be installed, the wind penetration level may increase but the value of objective function will decrease. To demonstrate the result, we further remove the investment term in the objective function (12a). The planning model indicates that seven TCSCs and seven PSTs will be installed. The wind penetration level increases to 42.53%, which is the highest penetration level that the system can achieve with series FACTS devices. The total cost, i.e., the overall sum of investment cost, wind curtailment cost and load shedding cost, is $197.61 M.

TABLE II
IEEE 118-BUS SYSTEM RESULTS FOR DIFFERENT CASES

| | Wind curtailment (×10⁶ MWh) | | | Load shedding (×10³ MWh) | TCSC locations | Investment on TCSC (M $) | PST locations | Investment on PST (M $) | Objective value (M $) | Wind penetration (%) | Time (s) |
|---|---|---|---|---|---|---|---|---|---|---|---|
| | 5 | 26 | 91 | | | | | | | | |
| $N_V=0$ $N_P=0$ | 1.1317 | 2.6570 | 1.7735 | 0.0056 | - | - | - | - | 278.1370 | 33.4056 | 1.7849 |
| $N_V=1$ $N_P=0$ | 1.1372 | 2.6556 | 1.0087 | 0 | 90-91 | 2.0656 | - | - | 242.1403 | 36.0675 | 16.4250 |
| $N_V=0$ $N_P=1$ | 1.1390 | 2.6552 | 0.7367 | 0 | - | - | 89-91 | 3.8111 | 230.3556 | 37.0144 | 2.7373 |
| $N_V=1$ $N_P=1$ | 1.0974 | 2.6986 | 0.4344 | 0 | 90-91 | 2.0656 | 89-91 | 3.8111 | 217.3941 | 38.0661 | 17.5054 |
| $N_V=2$ $N_P=1$ | 0.9938 | 2.6599 | 0.4345 | 0 | 90-91 19-34 | 3.4084 | 89-91 | 3.8111 | 211.6300 | 38.5635 | 53.3072 |
| $N_V=1$ $N_P=2$ | 0.8993 | 2.2832 | 0.7373 | 0 | 26-30 | 2.7670 | 89-91 30-38 | 7.6222 | 206.3771 | 39.1530 | 42.2448 |
| $N_V=2$ $N_P=2$ | 0.8993 | 2.2834 | 0.4352 | 0 | 90-91 26-30 | 4.8326 | 89-91 30-38 | 7.6222 | 193.3522 | 40.2091 | 20.5603 |

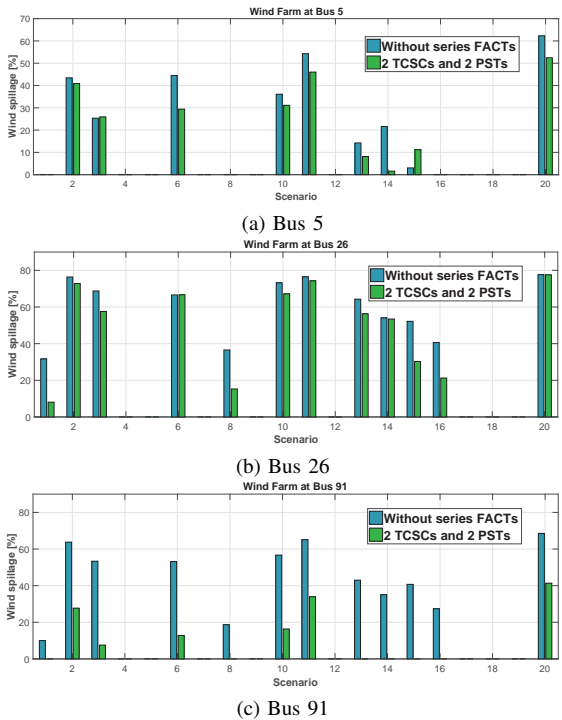

Fig. 5. Wind power spillage for each scenario.

Fig. 5 depicts the wind power curtailment ratio for the three wind farms in each scenario. It can be seen that the installation of four series FACTS reduces the curtailment ratios for most of the scenarios. The wind farm at bus 91 has the largest reductions. The number of scenarios in which the wind spillage occurs is 12 when no series FACTS is installed. With four series FACTS, only six scenarios involve the wind curtailment. Note that in scenario 15, the installation of four series FACTS decreases the use of wind power at bus 5 from 789.39 MW to 722.57 MW. Nevertheless, the wind power penetration for the other two wind farms both increase and the total wind power usage for this scenario increases from 1613.05 MW to 2022.66 MW.

In Table III, the series FACTS placement strategies provided by the proposed bilevel model and the conventional single-level model are compared for two cases. The single-level model minimizes the total cost, i.e., the summation of the investment cost and the operation cost, which is constrained by the budget and operation constraints. From Table III, it can be seen that the two approaches provide different locations for TCSC and PST. In addition, our proposed approach gives better planning results in terms of wind curtailment ratio. When $N_V=2$ and $N_P=2$, the ratio of wind curtailment is 27.16% by using the single-level model and this value decreases to 23.95% with the proposed bilevel model.

TABLE III
IEEE 118-BUS SYSTEM RESULTS COMPARISON

| | $N_V=1, N_P=1$ | | $N_V=2, N_P=2$ | |
|---|---|---|---|---|
| | Proposed model | Single-level model | Proposed model | Single-level model |
| TCSC locations | 90-91 | 19-34 | 90-91 26-30 | 26-30 19-34 |
| PST locations | 89-91 | 89-91 | 89-91 30-38 | 89-91 89-92 |
| Wind curtailment (%) | 27.9999 | 29.0548 | 23.9465 | 27.1558 |

### B. Computational Issues

All the simulations are conducted on a personal laptop with an Inter(R) Core(TM) i5-6300U CPU @ 2.40GHz and 8.00 GB of RAM. The complete model is implemented in MATLAB toolbox YALMIP [44] and solved by CPLEX [45]. The tolerance for the decomposition algorithm is set to 0.1%, the time limit for the master problem in one iteration is three hours.

To compare the computational performance, we consider four cases with respect to $N_V = 2$ and $N_P = 2$, the case description is given below:

C1: *Shift factor* structure with binary variable reduction in the lower level problems.

C2: $B\theta$ formulation with binary variable reduction in the lower level problems.

C3: *Shift factor* structure without binary variable reduction in the lower level problems.

C4: $B\theta$ formulation without binary variable reduction in the lower level problems.

Table IV provides the computational comparison results for different cases in terms of model size and computational time. The following observations can be draw:

- The planning results for the first three cases are exactly the same.
- The model size of *shift factor* structure dramatically decreases as compared to the $B\theta$ formulation because in *shift factor* structure: 1) bus angle variables vanish; 2) there is only one power balance equation; 3) The number of monitored transmission lines is 43 instead of 185. The *shift factor* formulation is beneficial to the algorithm since less new variables and constraints will be added to the master problem during the iterative process.
- Most computational time is spent on the master problem.
- The computational bottleneck for the algorithm is the number of binary variables in the lower level problems. Using *shift factor* formulation is capable of further reducing the computational time. In case 3, the *shift factor* formulation is able to find the optimal planning results in about half an hour even without the binary variable reduction strategy. However, the $B\theta$ formulation in case 4 fails to do so within the given time limit.

## VI. CONCLUSION

This paper presents a bilevel model to co-optimize the locations of VSR and PST considering high penetration of wind power. The proposed planning model seeks to identify the investment decisions on series FACTS within a market environment. To capture the intermittent nature of wind power, we consider a collection of lower problems to represent the market clearing under different load-wind scenarios. The resulting model is a stochastic MIP bilevel model with binary variables in both levels. A customized reformulation and decomposition algorithm is implemented to solve this challenging model. In addition, we compare the computational performance of $B\theta$ and *shift factor* formulation for the series FACTS allocation problem. The numerical results based on IEEE 118-bus system illustrate the significant benefits of series FACTS in the wind power integration. Also, the *shift factor* structure outperforms the $B\theta$ formulation in terms of computational speed due to its reduced model size.

TABLE IV
COMPUTATIONAL COMPARISON FOR DIFFERENT CASES

| | C1 | C2 | C3 | C4 |
|---|---|---|---|---|
| # binary in UL | 20 | 20 | 20 | 20 |
| # binary in LL | 40 | 40 | 200 | 200 |
| # continuous variables in LL | 4320 | 9520 | 4320 | 9520 |
| # inequality constraints in UL | 8 | 8 | 8 | 8 |
| # inequality constraints in LL | 8560 | 21640 | 8560 | 21640 |
| # equality constraints in LL | 1320 | 6520 | 1320 | 6520 |
| MP time (s) | 18.72 (s) | 44.72 (s) | 1883.56 (s) | 2.56 (h) |
| SP time (s) | 1.84 | 1.38 | 4.83 | 4.11 |
| Total time | 20.56 (s) | 46.20 (s) | 1888.39 (s) | 2.56 (h) |
| # iter. | 2 | 2 | 5 | 4 |
| Objective (M $) | 193.35 | 193.35 | 193.35 | 205.85 |

[1] The model size is based on the "initial" model, i.e., without any new variables and constraints in step 5).

[2] For case 4, the algorithm terminates in iteration 5 because the computational time of MP exceeds the time limit. The results reported in the table are for the first 4 iterations.

[3] Reference [13] and [27] provide the $B\theta$ formulation regarding the line flow with PST and VSR respectively.

Note that the proposed model adopts DC power flow, which ignores power loss and reactive power. Therefore, the proposed model is suitable for preparatory power network design. For more detailed voltage and angular stability limits analysis, the obtained FACTS locations can be further evaluated by using a full AC power flow model.

In Section IV-B2, a threshold of 60% for the line loading is used to decide which lines are monitored in the planning model. One direction of our future work is to derive an adaptive method which can dynamically adjust the threshold. In addition, reference [28] presents an alternative reformulation approach for the power flow on line with VSR. It would be interesting to investigate the reformulation in our proposed bilevel planning model.

## APPENDIX A
## CALCULATION OF THE LIMITS OF $\Delta b_k^V$

To derive the limit of $\Delta b_k^V$, we assume that the compensation range of VSR, i.e., $x_k^V$, is given as follows

$$x_{k,V}^{\min} \leq x_k^V \leq x_{k,V}^{\max} \qquad (28)$$

The range of $b_k^V$ in (3) can be computed as

$$-\frac{x_{k,V}^{\max}}{x_k(x_k + x_{k,V}^{\max})} \leq b_k^V \leq -\frac{x_{k,V}^{\min}}{x_k(x_k + x_{k,V}^{\min})} \qquad (29)$$

Thus, the maximum and minimum value of $\Delta b_k^V$ can be expressed as given in (30).

$$-\frac{x_{k,V}^{\max}}{x_k + x_{k,V}^{\max}} \leq \Delta b_k^V \leq -\frac{x_{k,V}^{\min}}{x_k + x_{k,V}^{\min}} \qquad (30)$$

## APPENDIX B
### $B\theta$ AND *shift factor* FORMULATION FOR DCOPF

Based on $B\theta$ formulation, the DCOPF for one snapshot can be expressed as

$$\min_{\Xi_1} \quad \sum_{n \in \mathcal{G}} a_n^g P_n^g \tag{31a}$$

subject to

$$P_k = \frac{\theta_i - \theta_j}{x_k}, \ \forall k \tag{31b}$$

$$\sum_{n \in \mathcal{G}_i} P_n^g - \sum_{m \in \mathcal{D}_i} P_m^d = \sum_{k \in f(i)} P_k - \sum_{k \in t(i)} P_k, \ \forall i \tag{31c}$$

$$P_n^{g,\min} \le P_n^g \le P_n^{g,\max}, \ \forall n \tag{31d}$$

$$-S_k^{\max} \le P_k \le S_k^{\max}, \ \forall k \tag{31e}$$

$$\theta_i = 0, \ \forall i \in \mathcal{B}_{ref} \tag{31f}$$

where $i, j$ denote the from and to end of branch $k$. The optimization variables of (31) comprises the elements in set $\Xi_1 = \{P_n^g, \theta_i, P_k\}$.

As mentioned in Section II, the PTDF is defined as the sensitivity of the power flow on line $k$ with respect to the power injection at bus $i$. Thus, an equivalent formulation of DCOPF based on *shift factor* structure can be obtained as follows

$$\min_{\Xi_2} \quad \sum_{n \in \mathcal{G}} a_n^g P_n^g \tag{32a}$$

subject to

$$P_k = \boldsymbol{H}(k,i)(\sum_{n \in \mathcal{G}_i} P_n^g - \sum_{m \in \mathcal{D}_i} P_m^d), \ \forall k, \forall i \tag{32b}$$

$$\sum_{n \in \mathcal{G}} P_n^g - \sum_{m \in \mathcal{D}} P_m^d = 0, \ \forall n, \forall m \tag{32c}$$

$$P_n^{g,\min} \le P_n^g \le P_n^{g,\max}, \ \forall n \tag{32d}$$

$$-S_k^{\max} \le P_k \le S_k^{\max}, \ \forall k \tag{32e}$$

The optimization variables of (32) include the elements in set $\Xi_2 = \{P_n^g, P_k\}$.

Assume the power network has $n_b$ buses, $n_l$ branches and $n_g$ generators, the model size of the above two formulations are compared in Table V.

### TABLE V
#### MODEL SIZE COMPARISON

| | $B\theta$ | Shift Factor |
|---|---|---|
| # variables | $n_b + n_l + n_g$ | $n_g + n_l$ |
| # equality constraints | $n_l + n_b + 1$ | $n_l + 1$ |
| # inequality constraints | $2n_l + 2n_g$ | $2n_l + 2n_g$ |

As can be observed from Table V, the *shift factor* formulation contains less variables and constraints since the bus angle variables $\theta$ are removed and there is only one power balance equation. Moreover, if a subset of transmission lines are monitored, the number of variables and constraints in *shift factor* formulation can be further decreased. Hence, the *shift factor* formulation has better scalability than the $B\theta$ formulation.

## REFERENCES

[1] (2015) Wind Vision: A New Era for Wind Power in the United States. [Online]. Available: http://www.energy.gov/sites/prod/files/WindVision_Report_final.pdf

[2] M. Milligan, K. Porter, E. DeMeo, P. Denholm, H. Holttinen, B. Kirby, N. Miller, A. Mills, M. O'Malley, M. Schuerger, and L. Söder, "Wind power myths debunked," *IEEE Power Energy Mag.*, vol. 7, no. 6, pp. 89–99, 2009.

[3] A. Nasri, A. J. Conejo, S. J. Kazempour, and M. Ghandhari, "Minimzing wind power spillage using an OPF with FACTS devices," *IEEE Trans. Power Syst.*, vol. 29, no. 5, pp. 2150–2159, Sep. 2014.

[4] S. Gerbex, R. Cherkaoui, and A. J. Germond, "Optimal location of multi-type FACTS devices in a power system by means of genetic algorithm," *IEEE Trans. Power Syst.*, vol. 16, no. 3, pp. 537–543, Aug. 2001.

[5] X. Zhang, C. Rehtanz, and B. Pal, *Flexible AC Transmission Systems: Modelling and Control.* London, UK: Springer, 1999.

[6] A. Dimitrovski, Z. Li, and B. Ozpineci, "Magnetic amplifier-based power flow controller," *IEEE Trans. Power Del.*, vol. 30, no. 4, pp. 1708–1714, Aug. 2015.

[7] Green Electricity Network Integration. [Online]. Available: http://arpa-e.energy.gov/?q=arpa-e-programs/geni

[8] M. Saravanan, S. M. R. Slochanal, R. Venkatesh, and J. P. S. Abrahamh, "Application of PSO technique for optimal location of FACTS devices considering system loadability and cost of installation," *Elect. Power Syst. Res.*, vol. 77, no. 3, pp. 276–283, Mar. 2007.

[9] R. S. Wibowo, N. Yorino, M. Eghbal, Y. Zoka, and Y. Sasaki, "FACTS devices allocation with control coordination considering congestion relief and voltage stability," *IEEE Trans. Power Syst.*, vol. 26, no. 4, pp. 2302–2310, Nov. 2011.

[10] A. R. Jordehi and J. Jasni, "Heuristic methods for solution of FACTS optimization problem in power systems," in *Proc. 2011 IEEE Student Conf. on Research and Development (SCOReD)*, Cyberjaya, Malaysia, Dec. 19–20, 2011, pp. 30–35.

[11] T. Orfanogianni and R. Bacher, "Steady-state optimization in power systems with series FACTS devices," *IEEE Trans. Power Syst.*, vol. 18, no. 1, pp. 19–26, Feb. 2003.

[12] Y. Lu and A. Abur, "Static security enhancement via optimal utilization of thyristor-controlled series capacitors," *IEEE Trans. Power Syst.*, vol. 17, no. 2, pp. 324–329, May 2002.

[13] F. G. M. Lima, F. D. Galiana, I. Kockar, and J. Munoz, "Phase shifter placement in large-scale system via mixed integer linear programming," *IEEE Trans. Power Syst.*, vol. 18, no. 3, pp. 1029–1034, Aug. 2003.

[14] G. Yang, G. Hovland, R. Majumder, and Z. Dong, "TCSC allocation based on line flow based equations via mixed-integer programming," *IEEE Trans. Power Syst.*, vol. 22, no. 4, pp. 2262–2269, Nov. 2007.

[15] X. Zhang, K. Tomsovic, and A. Dimitrovski, "Optimal investment on series FACTS device considering contingencies," in *Proc. North Amer. Power Symp., 2016*, Denver, CO, USA, Sep. 18–20, 2016, pp. 1–6.

[16] O. Ziaee and F. Choobineh, "Optimal location-allocation of TCSCs and transmission switch placement under high penetration of wind power," *IEEE Trans. Power Syst.*, vol. 32, no. 4, pp. 3006–3014, Jul. 2017.

[17] L. P. Garcés, A. J. Conejo, R. García-Bertrand, and R. Romero, "A bilevel approach to transmission expansion planning within a market environment," *IEEE Trans. Power Syst.*, vol. 24, no. 3, pp. 1513–1522, Aug. 2009.

[18] F. Ugranli, E. Karatepe, and A. H. Nielsen, "MILP approach for bilevel transmission and reactive power planning considering wind curtailment," *IEEE Trans. Power Syst.*, vol. 32, no. 1, pp. 652–661, Jan. 2017.

[19] L. Baringo and A. Conejo, "Wind power investment within a market environment," *App. Energy*, vol. 88, no. 9, pp. 3239–3247, Sep. 2011.

[20] S. J. Kazempour, A. J. Conejo, and C. Ruiz, "Strategic generation investment using a complementarity approach," *IEEE Trans. Power Syst.*, vol. 26, no. 2, pp. 940–948, May 2011.

[21] L. Baringo and A. J. Conejo, "Transmission and wind power investment," *IEEE Trans. Power Syst.*, vol. 27, no. 2, pp. 885–893, May 2012.

[22] T. Ding, R. Bo, Z. Bie, and X. Wang, "Optimal selection of phase shifting transformer adjustment in optimal power flow," *IEEE Trans. Power Syst.*, vol. 32, no. 3, pp. 2464–2465, May 2017.

[23] M. Rahmani, A. Kargarian, and G. Hug, "Comprehensive power transfer distribution factor model for large-scale transmission expansion planning," *IET Gener., Transm., Distrib.*, vol. 10, no. 12, pp. 2981–2989, 2016.

[24] M. Sahraei-Ardakani and K. W. Hedman, "Computationally efficient adjustment of FACTS set points in dc optimal power flow with shift factor structure," *IEEE Trans. Power Syst.*, vol. 32, no. 3, pp. 1733–1740, May 2017.

[25] B. Zeng and Y. An, "Solving bilevel mixed integer program by reformulations and decomposition," *Optimization online*, pp. 1–34, 2014.

[26] T. T. Lie and W. Deng, "Optimal flexible AC transmission systems (FACTS) devices allocation," *Electr. Power Energy Syst.*, vol. 19, no. 2, pp. 125–134, Feb. 1997.

[27] X. Zhang, K. Tomsovic, and A. Dimitrovski, "Security constrained multi-stage transmission expansion planning considering a continuously variable series reactor," *IEEE Trans. Power Syst.*, vol. 32, no. 6, pp. 4442–4450, Nov. 2017.

[28] O. Ziaee, O. Alizadeh-Mousavi, and F. F. Choobineh, "Co-optimization of transmission expansion planning and TCSC placement considering the correlation between wind and demand scenarios," *IEEE Trans. Power Syst.*, vol. 33, no. 1, pp. 206–215, Jan. 2018.

[29] B. Zeng and L. Zhao, "Solving two-stage robust optimization problems using a column-and-constraint generation method," *Oper. Res. Lett.*, vol. 41, no. 5, pp. 457–461, Sep. 2013.

[30] W. Yuan, J. Wang, F. Qiu, C. Chen, C. Kang, and B. Zeng, "Robust optimization-based resilient distribution network planning against natural disasters," *IEEE Trans Smart Grid*, vol. 7, no. 6, pp. 2817–2826, Nov. 2016.

[31] H. Haghighat and B. Zeng, "Distribution system reconfiguration under uncertain load and renewable generation," *IEEE Trans. Power Syst.*, vol. 31, no. 4, pp. 2666–2675, 2016.

[32] IPOPT. [Online]. Available: https://projects.coin-or.org/Ipopt

[33] M. S. Ardakani and K. W. Hedman, "A fast LP approach for enhanced utilization of variable impedance based FACTS devices," *IEEE Trans. Power Syst.*, vol. 31, no. 3, pp. 2204–2213, May 2016.

[34] PLEXOS wiki. [Online]. Available: https://wiki.energyexemplar.com/

[35] S. A. Blumsack, "Network topologies and transmission investment under electric-industry restructuring," Ph.D. dissertation, Carnegie Mellon Univ., Pittsburgh, PA, May 2006. [Online]. Available: http://www.personal.psu.edu/sab51/Blumsack_Dissertation.pdf

[36] L. Ippolito and P. Siano, "Selection of optimal number and location of thyristor-controlled phase shifters using genetic based algorithms," *IEE Proc. Gener., Transm., Distrib.*, vol. 151, no. 5, pp. 630–637, Sep. 2004.

[37] O. Ziaee and F. Choobineh, "Stochastic location-allocation of TCSC devices on a power system with large scale wind generation," in *Proc. IEEE Power Eng. Soc. Gen. Meeting*, Boston, MA, USA, Jul. 17–21, 2016, pp. 1–5.

[38] Y. Sang and M. S. Ardakani, "The interdependence between transmission switching and variable-impedance series facts devices," *IEEE Trans. Power Syst.*, vol. 33, no. 3, pp. 2792–2803, Sep. 2017.

[39] P. K. Tiwari and Y. R. Sood, "An efficient approach for optimal allocation and parameters determination of TCSC with investment cost recovery under competitive power market," *IEEE Trans. Power Syst.*, vol. 28, no. 3, pp. 2475–2484, Aug. 2013.

[40] Energy, Load, and Demand Reports. [Online]. Available: https://www.iso-ne.com/isoexpress/web/reports/load-and-demand/-/tree/dmnd-rt-hourly-sys

[41] Renewables.ninja. [Online]. Available: https://www.renewables.ninja/

[42] L. Baringo and A. Conejo, "Correlated wind-power production and electric load scenarios for investment decisions," *App. Energy*, vol. 101, pp. 475–482, Jan. 2013.

[43] J. C. Villumsen, G. Bronmo, and A. B. Philpott, "Line capacity expansion and transmission switching in power systems with large-scale wind power," *IEEE Trans. Power Syst.*, vol. 28, no. 2, pp. 731–739, May 2013.

[44] J. Löfberg, "YALMIP: A toolbox for modeling and optimization in MATLAB," in *Proc. CACSD Conf.*, Taipei, Taiwan, Jul. 2004, pp. 284–289.

[45] (2014) IBM ILOG CPLEX V 12.6. [Online]. Available: http://www.ibm.com/software/commerce/optimization/cplex-optimizer/

**Xiaohu Zhang** (S'12, M'17) received the B.S. degree in electrical engineering from Huazhong University of Science and Technology, Wuhan, China, in 2009, the M.S. degree in electrical engineering from Royal Institute of Technology, Stockholm, Sweden, in 2011, and the Ph.D. degree in electrical engineering at The University of Tennessee, Knoxville, in 2017.

Currently, he works as a power system engineer at GEIRI North America, San Jose, CA, USA. His research interests are power system operation, planning and stability analysis.

**Di Shi** (M'12, SM'17) received the Ph.D. degree in electrical engineering from Arizona State University, Tempe, AZ, USA, in 2012. He currently leads the PMU & System Analytics Group at GEIRI North America, San Jose, CA, USA. Prior to that, he was a researcher at NEC Laboratories America, Cupertino, CA, and Electric Power Research Institute (EPRI), Palo Alto, CA. He served as Senior/Principal Consultant for eMIT and RM Energy Marketing between 2012-2016. He has published over 70 journal and conference papers and holds 14 US patents/patent applications. He received the IEEE PES General Meeting Best Paper Award in 2017. One Energy Management and Control (EMC) technology he developed has been commercialized in 2014 into product that helps customers achieve significant energy savings. He is an Editor of IEEE Transactions on Smart Grid.

**Zhiwei Wang** received the B.S. and M.S. degrees in electrical engineering from Southeast University, Nanjing, China, in 1988 and 1991, respectively. He is President of GEIRI North America, San Jose, CA, USA. His research interests include power system operation and control, relay protection, power system planning, and WAMS.

**Bo Zeng** (M'11) received the Ph.D. degree in industrial engineering from Purdue University, West Lafayette, IN, USA, with emphasis on operations research. He is currently an assistant professor in the Department of Industrial Engineering and the Department of Electrical and Computer Engineering, University of Pittsburgh, Pittsburgh, PA, USA. His research interests include polyhedral theory and algorithm development for stochastic, robust, and multilevel mixed-integer programs, with applications in power, logistics and cyber-physical systems. He is a member of IISE, INFORMS, and IEEE.

**Xinan Wang** (S'15) received the B.S. degree in electrical engineering from Northwestern Polytechnical University, Xian, China, in 2013, and the M.S. degree in electrical engineering from Arizona State University, Tempe, AZ, USA, in 2016. During Aug. 2016-May. 2017, he worked as a research assistant in the Advanced Power System Analytics Group at GEIRI North America, Santa Clara, CA, USA. He is currently pursuing his PhD degree in electrical engineering at Southern Methodist University, Dallas, TX, USA. His research interests include WAMS related application in power system, data driven load monitoring and renewable energy integration.

**Kevin Tomsovic** (F'07) received the BS from Michigan Tech. University, Houghton, in 1982, and the MS and Ph.D. degrees from University of Washington, Seattle, in 1984 and 1987, respectively, all in Electrical Engineering. He is currently the CTI Professor with the Department of Electrical Engineering and Computer Science, University of Tennessee, Knoxville, TN, USA, where he directs the NSF/DOE ERC, Center for Ultra-Wide-Area Resilient Electric Energy Transmission Networks (CURENT), and previously served as the Electrical Engineering and Computer Science Department Head from 2008 to 2013. He was on the faculty of Washington State University, Pullman, WA, USA, from 1992 to 2008. He held the Advanced Technology for Electrical Energy Chair at National Kumamoto University, Kumamoto, Japan, from 1999 to 2000, and was the NSF Program Director with the Electrical and Communications Systems Division of the Engineering Directorate from 2004 to 2006. He also held positions at National Cheng Kung University and National Sun Yat Sen University in Taiwan from 1988-1991 and the Royal Instituted of Technology in Sweden from 1991-1992. He is a Fellow of the IEEE.

**Yanming Jin** received the Ph.D. degree from Guanghua School of Management, Peking University, Beijing, China, in 2007. She is with State Grid Energy Research Institute, Beijing, China. Her research interests includes energy related strategy, planning and environmental policy.